\documentclass[12pt]{article}
\usepackage{amsmath, amsthm, amssymb, stmaryrd}
\usepackage[all]{xy}
\usepackage{fullpage}
\usepackage{titlesec}

\titleformat{\section}{\normalsize\bfseries}{\thesection}{1em}{}
\titleformat{\subsection}{\normalsize\bfseries}{\thesubsection}{1em}{}

\theoremstyle{plain}

\newtheorem{PropSub}[subsection]{Proposition}
\newtheorem{LemSub}[subsection]{Lemma}

\newtheorem{ThmSub}[subsection]{Theorem}

\newtheorem{PropSubSub}[subsubsection]{Proposition}

\newtheorem{ThmSubSub}[subsubsection]{Theorem}

\theoremstyle{definition}

\newtheorem{DefSub}[subsection]{Definition}

\newtheorem{RemSubSub}[subsubsection]{Remark}
\newtheorem{ExaSubSub}[subsubsection]{Example}

\newcommand*{\emptybox}{\leavevmode\hbox{}}

\DeclareMathOperator*{\col}{colim}

\newcommand{\cl}[1]{\ensuremath{\displaystyle \col_{#1}}}

\newcommand{\h}[1]{\ensuremath{\widehat{#1}}}

\newcommand{\A}{\ensuremath{\mathcal{A}}}

\newcommand{\C}{\ensuremath{\mathcal{C}}}
\newcommand{\D}{\ensuremath{\mathcal{D}}}
\newcommand{\E}{\ensuremath{\mathcal{E}}}
\newcommand{\F}{\ensuremath{\mathcal{F}}}
\newcommand{\G}{\ensuremath{\mathcal{G}}}

\newcommand{\X}{\ensuremath{\mathcal{X}}}

\newcommand{\CAT}{\ensuremath{\operatorname{\textnormal{\textbf{CAT}}}}}
\newcommand{\CCAATT}{\ensuremath{\mathfrak{CAT}}}
\newcommand{\Op}[1]{\ensuremath{#1^{\textnormal{op}}}}

\newcommand{\Prof}{\ensuremath{\operatorname{\textnormal{\textbf{Prof}}}}}

\newcommand{\pt}{\ensuremath{\operatorname{\textnormal{\textsf{pt}}}}}

\newcommand{\Dn}{\ensuremath{\operatorname{\textnormal{\textsf{Dn}}}}}
\newcommand{\Idl}{\ensuremath{\operatorname{\textnormal{\textsf{Idl}}}}}
\newcommand{\Ind}{\ensuremath{\operatorname{\textnormal{\textsf{Ind}}}}}

\newcommand{\y}{\ensuremath{\textnormal{\textsf{y}}}}
\newcommand{\m}{\ensuremath{\textnormal{\textsf{m}}}}

\newcommand{\twoheaddownarrow}{\mathord{\mbox{\makebox[0pt][l]{$\downarrow$}\raisebox{0.4ex}{$\downarrow$}}}}

\newcommand{\threeheaddownarrow}{\mathord{\mbox{\makebox[0pt][l]{$\downarrow$}\raisebox{0.4ex}{$\twoheaddownarrow$}}}}

\newcommand{\Set}{\ensuremath{\operatorname{\textnormal{\textbf{Set}}}}}

\def\mrelarrow{\mathop{\rlap{\hspace{0.5ex}$\circ$}{\rightarrow}}\nolimits}

\newcommand{\op}{\ensuremath{\textnormal{op}}}
\newcommand{\coop}{\ensuremath{\textnormal{coop}}}

\newcommand{\lt}{\leqslant}

\begin{document}

\author{\normalsize  Rory B.B. Lucyshyn-Wright\thanks{Partial financial assistance by the Ontario Graduate Scholarship program is gratefully acknowledged.}\let\thefootnote\relax\let\thefootnote\relax
\\
\small York University, 4700 Keele St., Toronto, ON, Canada M3J 1P3}

\title{\large \textbf{Totally distributive toposes}}

\date{}

\maketitle

\abstract{A locally small category $\E$ is \textit{totally distributive} (as defined by Rosebrugh-Wood) if there exists a string of adjoint functors $t \dashv c \dashv \y$, where $\y:\E \rightarrow \h{\E}$ is the Yoneda embedding.  Saying that $\E$ is \textit{lex totally distributive} if, moreover, the left adjoint $t$ preserves finite limits, we show that the lex totally distributive categories with a small set of generators are exactly the \textit{injective Grothendieck toposes}, studied by Johnstone and Joyal.  We characterize the totally distributive categories with a small set of generators as exactly the \textit{essential} subtoposes of presheaf toposes, studied by Kelly-Lawvere and Kennett-Riehl-Roy-Zaks.
}

\section{Introduction} \label{sec:intro}

The aim of this paper is to establish certain connections between the work of Marmolejo, Rosebrugh, and Wood \cite{RWsets,MRW} on \textit{totally distributive categories} and two other bodies of work on distinct topics:  Firstly, that of Johnstone and Joyal \cite{Jinj,JJ} on \textit{injective toposes} and \textit{continuous categories}, and secondly, that of Kelly-Lawvere \cite{KL} and Kennett-Riehl-Roy-Zaks \cite{KRRZ} on \textit{essential localizations} and \textit{essential subtoposes}.  One of our observations, \ref{thm:main} (2), when taken together with a theorem of Kelly-Lawvere which we recall in \ref{rem:KL_thm}, yields a concrete combinatorial description of all totally distributive categories with a small set of generators.

We adopt the foundational conventions of $\cite{Jele}$ (and \cite{Jinj,JJ}), since our only use of the stronger foundational assumptions of $\cite{SW,St,W,RWsets,MRW}$ is made in finally deducing our main results (\ref{thm:main}) as strengthened variants of propositions which precede them.  We let $\CCAATT$ represent the meta-2-category of categories, functors, and natural transformations (see \cite{Jele}, 1.1.1), and we let $\CAT$ be its full sub-(meta)-2-category consisting of locally small categories.

\renewcommand*\theparagraph{\thesection.\arabic{subsection}}
\addtocounter{subsection}{1}

\paragraph{Completely distributive lattices, totally distributive categories.} \label{sec:intro_cd_lats_td_cats}

A poset $\E$ is a \textit{constructively completely distributive lattice} \cite{FW}, or \textit{ccd lattice}, if there exist adjunctions
$$
\xymatrix {
\E \ar@/_1.2pc/[rr]_\threeheaddownarrow^{\top} \ar@/^1.2pc/[rr]^\downarrow_{\top} & & \Dn(\E) \ar[ll]|\vee
}
$$
where $\Dn(\E)$ is the poset of down-closed subsets of $\E$, ordered by inclusion, and $\downarrow:\E \rightarrow \Dn(\E)$ is the embedding given by $v \mapsto \;\downarrow v := \{u \in \E \mid u \lt v\}$.  The existence of the left adjoint $\vee$ of $\downarrow$ is equivalent to the cocompleteness of $\E$, i.e. the condition that $\E$ be a complete lattice, and if such a map $\vee$ exists, it necessarily sends each down-closed subset to its join in $\E$.  In the presence of the axiom of choice, a poset is a ccd lattice iff it is a completely distributive lattice in the usual sense \cite{FW}.

Rosebrugh and Wood \cite{RWsets} have defined an analogue of this notion for arbitrary categories rather than just posets\footnote{Marmolejo, Rosebrugh, and Wood \cite{MRW} have also studied an apparently distinct analogue --- the notion of \textit{completely distributive category}.}.  A locally small category $\E$ is \textit{totally distributive} if there exist adjunctions
$$
\xymatrix {
\E \ar@/_1.1pc/[rr]_t^{\top} \ar@/^1.1pc/[rr]^\y_{\top} & & \h{\E} \ar[ll]|c
}
$$
where $\h{\E}$ is the presheaf category $[\E^\op, \Set]$ and $\y$ is the Yoneda embedding, given by $v \mapsto \h{v} := \E(-,v)$.  We say that a totally distributive category $\E$ is \textit{lex totally distributive} if the associated functor $t:\E \rightarrow \h{\E}$ preserves finite limits.

The existence of the left adjoint $c$ of $\y$ is the requirement that $\E$ be \textit{total} \cite{SW}, or \textit{totally cocomplete}.  This left adjoint $c$ of $\y$ is characterized by the property that
\begin{equation}\label{eq:col_of_presh}cE \cong \cl{\h{u} \rightarrow E} u = \col((\E \downarrow E) \rightarrow \E) \cong \int^{u \in \E}Eu \cdot u\end{equation}
naturally in $E \in \h{\E}$, so that totality is equivalent to the existence of a colimit in $\E$ of the (possibly large) canonical diagram of each presheaf $E$ on $\E$.

Note that any totally distributive category $\E$ is in particular \textit{lex total}, meaning that $\E$ is total and the associated functor $c:\h{\E} \rightarrow \E$ preserves finite limits.  Wood \cite{W} attributes to Walters the theorem that those lex total categories with a small set of generators are exactly the Grothendieck toposes; the paper \cite{St} of Street includes a proof of this result.

\addtocounter{subsection}{1}
\paragraph{Continuous dcpos, continuous categories.}

A poset $\X$ is a \textit{continuous dcpo} if there exist adjunctions
$$
\xymatrix {
\X \ar@/_1.2pc/[rr]_\twoheaddownarrow^{\top} \ar@/^1.2pc/[rr]^\downarrow_{\top} & & \Idl(\X) \ar[ll]|\vee
}
$$
where $\Idl(\X)$ is the poset of ideals of $\X$ (i.e. upward-directed down-closed subsets of $\X$), ordered by inclusion, and $\downarrow:\X \rightarrow \Idl(\X)$ is the embedding given by $y \mapsto \;\downarrow y := \{x \in \X \mid x \lt y\}$.  The existence of the left adjoint $\vee$ of $\downarrow$ is equivalent to the existence of all directed joins in $\X$, i.e. the condition that $\X$ be a \textit{dcpo}, or \textit{directed complete partial order}, and if such a map $\vee$ exists, it necessarily sends each ideal to its join in $\E$.

Johnstone and Joyal \cite{JJ} have defined a generalization of this notion to arbitrary categories, rather than just posets, as follows.  We say that a locally small category $\X$ is \textit{continuous} if there exist adjunctions
$$
\xymatrix {
\X \ar@/_1.1pc/[rr]_w^{\top} \ar@/^1.1pc/[rr]^\m_{\top} & & \Ind{\X} \ar[ll]|\col
},
$$
where $\Ind{\X}$ is the \textit{ind-completion} of $\X$, whose objects are all small filtered diagrams in $\X$, and $\m$ is the canonical full embedding sending each object $x \in \X$ to the diagram $1 \rightarrow \X$, indexed by the terminal category $1$, with constant value $x$.

The existence of the left adjoint $\col$ of $\m:\X \rightarrow \Ind{X}$ is equivalent to the requirement that $\X$ be equipped with colimits for all small filtered diagrams, and $\col$ necessarily sends each $D \in \Ind{\X}$ to a colimit of $D$ in $\X$.

\addtocounter{subsection}{1}
\paragraph{Stone duality for continuous dcpos.} \label{sec:stone_duality_for_cts_dcpos}

It was shown by Hoffmann \cite{H} and Lawson \cite{L} that the category  of continuous dcpos and directed-meet-preserving maps is dually equivalent to the category of completely distributive lattices and maps preserving finite meets and arbitrary joins.  
The category of continuous dcpos is isomorphic to the full subcategory of topological spaces consisting of continuous dcpos endowed with the \textit{Scott topology}, and the given dual equivalence of this category of spaces with the category of completely distributive lattices is a restriction of the dual equivalence between sober spaces and spatial frames (see, e.g., \cite{Jsto}), associating to a space its frame of open sets.

Subsequent work of Banaschewski \cite{B} entails that this dual equivalence restricts further to a dual equivalence between \textit{continuous lattices} (i.e. those continuous dcpos which are also complete lattices) and \textit{stably supercontinuous lattices}, also known as \textit{lex ccd lattices} \cite{MRW} or \textit{lex completely distributive lattices}, which are those ccd lattices for which the left adjoint $\threeheaddownarrow$ preserves finite meets.  Scott \cite{Sc} had shown earlier that the continuous lattices, when endowed with their Scott topologies, are exactly the \textit{injective} $T_0$ spaces.

\addtocounter{subsection}{1}
\paragraph{Continuous categories and injective toposes.}

Scott's isomorphism between injective $T_0$ spaces and continuous lattices \cite{Sc} has a topos-theoretic analogue, given by Johnstone-Joyal \cite{JJ}, which we now recall.

First let us record the following earlier result of Johnstone \cite{Jinj}:

\begin{ThmSubSub} \label{thm:injtop_as_retracts}
\emph{(Johnstone \cite{Jinj}).}
A Grothendieck topos $\E$ is injective (with respect to geometric inclusions) if and only if $\E$ is a retract, by geometric morphisms, of a presheaf topos $\h{\C}$ with $\C$ a small finitely-complete category.
\end{ThmSubSub}

We call such Grothendieck toposes \textit{injective toposes}.  A \textit{quasi-injective topos} \cite{JJ} is defined as a Grothendieck topos which is a retract, by geometric morphisms, of an arbitrary presheaf topos $\h{\C}$ (with $\C$ a small category).  A continuous category $\X$ is \textit{ind-small} if there exists a small \textit{ind-dense} subcategory $\A$ of $\X$, by which we mean a small, full, dense subcategory $\A$ of $\X$ for which each comma category $(\A \downarrow x)$, with $x \in \X$, is filtered\footnote{The term \textit{ind-small} was introduced not in \cite{JJ} but later in \cite{Jele}, where it is defined in terms of a different criterion, which, by 2.17 of \cite{JJ} and C4.2.18 of \cite{Jele}, is equivalent to the given condition, employed in \cite{JJ}.  Chapter C4 of \cite{Jele} includes an alternate exposition of much of the content of \cite{JJ}.}.

\begin{ThmSubSub} \label{thm:injcts}
\emph{(Johnstone-Joyal \cite{JJ}).}
\begin{enumerate}
\item There is an equivalence of 2-categories between the 2-category of quasi-injective toposes, with geometric morphisms, and the 2-category of ind-small continuous categories, with morphisms all filtered-colimit-preserving functors.  This equivalence sends a quasi-injective topos $\E$ to its category of points $\pt(\E)$.
\item This equivalence restricts to an equivalence between the full sub-2-categories of injective toposes and cocomplete ind-small continuous categories.
\end{enumerate}
\end{ThmSubSub}

\addtocounter{subsection}{1}
\paragraph{Totally distributive toposes.}\label{sec:qinjtop}

Having seen that Scott's isomorphism between injective $T_0$ spaces and continuous lattices has a topos-theoretic analogue relating injective toposes and cocomplete ind-small continuous categories, we are led to seek a topos-theoretic analogue of the dual equivalence between continuous lattices and lex completely distributive lattices.  We prove the following, where by a \textit{small dense generator} for a category $\E$ we mean a small dense full subcategory $\G$ of $\E$.  Recall that every Grothendieck topos has a small dense generator.

\begin{ThmSubSub} \label{thm:lex_td_sdg_equals_inj_top}
The lex totally distributive categories with a small dense generator are exactly the injective toposes.  Hence, the 2-category of cocomplete ind-small continuous categories (\ref{thm:injcts}) is equivalent to the the 2-category of lex totally distributive categories with a small dense generator (with geometric morphisms).
\end{ThmSubSub}

One may also ask whether there is a similar analogue of the broader dual equivalence between continuous dcpos and completely distributive lattices, and in this regard we provide a partial result, as follows:

\begin{PropSubSub} \label{thm:qinj_implies_td}
Every quasi-injective topos is totally distributive.
\end{PropSubSub}

In proving these theorems, we come upon a further result of independent interest.  An \textit{essential subtopos} of a topos $\F$ is a topos $\E$ for which there is a geometric inclusion $i:\E \rightarrow \F$ whose inverse image functor $i^*:\F \rightarrow \E$ has a left adjoint.

\begin{ThmSubSub} \label{thm:charn_of_td_sdg}
Those totally distributive categories having a small dense generator are exactly the essential subtoposes of presheaf toposes $\h{\C} = [\C^\op,\Set]$ (with $\C$ a small category).
\end{ThmSubSub}

\begin{RemSubSub} \label{rem:KL_thm}
It was shown by Kelly-Lawvere \cite{KL} that the essential subtoposes of a presheaf topos $\h{\C}$ correspond bijectively to \textit{idempotent ideals} of arrows in the small category $\C$.
\end{RemSubSub}

\begin{ExaSubSub}
The cases in which $\h{\C}$ is the topos $\h{\Delta}$ of \textit{simplicial sets}, the topos $\h{\mathbb{I}}$ of \textit{cubical sets}, or the topos $\h{\mathbb{G}}$ of \textit{reflexive globular sets} are of interest in homotopy theory and higher category theory.  It is shown in \cite{KRRZ} that the essential subtoposes of these toposes are classified by the dimensions $n \in \mathbb{N}$.  In general, the essential subtoposes of a topos $\F$ (or rather, their associated equivalent full replete subcategories of $\F$) form a complete lattice \cite{KL}.
\end{ExaSubSub}

\begin{RemSubSub} \label{rem:small_set_gens}
As noted in \ref{sec:intro_cd_lats_td_cats}, it was proved in \cite{St} that any lex total category $\E$ with a small set of generators is a Grothendieck topos.  Using this result, whose proof in \cite{St} appears to make use of the foundational assumption that there is a category of sets $S'$ such that both $\E$ and the category $\Set$ of small sets are categories internal to $S'$, we obtain the following corollaries to Theorems \ref{thm:lex_td_sdg_equals_inj_top} and \ref{thm:charn_of_td_sdg}:
\end{RemSubSub}

\begin{ThmSubSub} \label{thm:main}
\emptybox

\begin{enumerate}
\item Those lex totally distributive categories having a small set of generators are exactly the injective toposes.
\item Those totally distributive categories having a small set of generators are exactly the essential subtoposes of presheaf toposes $\h{\C} = [\C^\op,\Set]$ (with $\C$ small).
\end{enumerate}
\end{ThmSubSub}

\renewcommand*\theparagraph{\thesubsubsection.\arabic{paragraph}}

\section{Preliminaries on totally distributive categories}

It is shown in \cite{RWsets}, by means of a result of \cite{SW}, that every presheaf category $\h{\C}$ on a small category $\C$ is totally distributive.  In order to clearly establish this in the absence of the foundational assumptions of \cite{RWsets}, we give a self-contained elementary proof, by means of the following lemma (cf. Corollary 14 of \cite{SW}).  We prove also that if $\C$ is finitely complete, then $\h{\C}$ is lex totally distributive.

\begin{LemSub} \label{thm:presh_lemma}
Let $\C$ be a small category.  Then there is an adjunction
$$
\xymatrix {
\h{\C} \ar@/^1.2pc/[rr]^{\y_{\h{\C}}}_{\top} & & \h{\h{\C}} \ar[ll]|{\widehat{\y_\C}}
}\;,
$$
where $\y_\C:\C \rightarrow \h{\C}$ and $\y_{\h{\C}}:\h{\C} \rightarrow \h{\h{\C}}$ are the Yoneda embeddings.
\end{LemSub}
\begin{proof}
Each $\mathbb{C} \in \h{\h{\C}}$ is a coend $\mathbb{C} \cong \int^{C \in \h{\C}}\mathbb{C}(C) \cdot \h{C}$, and these isomorphisms are natural in $\mathbb{C}$.  Using this and the Yoneda Lemma, we obtain isomorphisms
$$(\widehat{\y_{\C}}(\mathbb{C}))(c) = \mathbb{C}(\h{c}) \cong \int^{C \in \h{\C}}\mathbb{C}(C) \times \h{C}(\h{c}) \cong \int^{C \in \h{\C}}\mathbb{C}(C) \times C(c)$$
natural in $\mathbb{C} \in \h{\h{\C}}$ and $c \in \C$.  Hence we have an isomorphism
$$\widehat{\y_{\C}}(\mathbb{C}) \cong \int^{C \in \h{\C}}\mathbb{C}(C) \cdot C$$
natural in $\mathbb{C} \in \h{\h{\C}}$, so with reference to \eqref{eq:col_of_presh}, $\widehat{\y_{\C}} \dashv \y_{\h{\C}}$.
\end{proof}

\begin{PropSub} \label{thm:preshtd}
Let $\C$ be a small category.  Then $\h{\C}$ is totally distributive.  Moreover, if $\C$ has finite limits, then $\h{\C}$ is lex totally distributive.
\end{PropSub}
\begin{proof}
We have an adjunction as in Lemma \ref{thm:presh_lemma}, and the left adjoint $\widehat{\y_\C}:\h{\h{\C}} \rightarrow \h{\C}$ has a further left adjoint $\exists_{\y_\C}:[\C^\op,\Set] \rightarrow [{\h{\C}\,}^\op,\Set]$, which is given by left Kan-extension along $\y_\C^\op:\C^\op \rightarrow {\h{\C}\,}^\op$.  Hence $\h{\C}$ is totally distributive.  If $\C$ has finite limits, then $\y_\C:\C \rightarrow \h{\C}$ is a cartesian functor between cartesian categories, and it follows that the associated functor $\exists_{\y_\C}$ is also cartesian.
\end{proof}
The following lemma, based on Lemma 3.5 of Marmolejo-Rosebrugh-Wood \cite{MRW}, provides a means of deducing that a category is totally distributive.  We have augmented the lemma slightly in order to handle lex totally distributive categories as well.

\begin{LemSub} \label{thm:mrwlemma}
Let $\D$ and $\E$ be locally small categories.  Suppose we are given adjunctions
$$
\xymatrix {
\D \ar@/_1.1pc/[rr]_q^{\top} \ar@/^1.1pc/[rr]^s_{\top} & & \E \ar[ll]|r
}
$$
with $q,s$ fully faithful and $\E$ totally distributive.  Then $\D$ is totally distributive.

Moreover, if $\E$ is lex totally distributive and $q$ preserves finite limits, then $\D$ is lex totally distributive.
\end{LemSub}
\begin{proof}
There is a 2-functor $\h{(-)} := \CAT((-)^\op,\Set) : \CAT^\coop \rightarrow \CCAATT$, where $\CAT^\coop$ is the (meta)-2-category gotten by reversing both the 1-cells and 2-cells in $\CAT$.  This 2-functor sends the adjunctions $q \dashv r \dashv s : \D \rightarrow \E$ in $\CAT$ to adjunctions $\h{q} \dashv \h{r} \dashv \h{s}$, so we have a diagram
$$
\xymatrix @+1pc {
\D \ar@/_1.0pc/[rr]_q^{\top} \ar@/^1.0pc/[rr]^s_{\top} \ar@/^1.0pc/[d]^{\y'} & & \E \ar[ll]|r \ar@{->}@/_1.0pc/[d]_t^{\dashv} \ar@{->}@/^1.0pc/[d]^{\y}_{\dashv}\\
\h{\D} \ar@{<-}@/_1.1pc/[rr]_{\h{q}}^{\top} \ar@{<-}@/^1.1pc/[rr]^{\h{s}}_{\top} & & \h{\E} \ar@{<-}[ll]|{\h{r}} \ar[u]|c
}
$$
where $\y'$ is the Yoneda embedding.  Observe that $\y' \cong \h{s} \cdot \y \cdot s$, since we have
$$(\h{s} \cdot \y \cdot s)(d) = \h{s}(\E(-,sd)) = \E(\Op{s}-,sd) \cong \D(-,d) = \y'(d) \;\;\;\;$$
naturally in $d \in \D$, as $s$ is fully faithful.  Therefore, letting $c' := r \cdot c \cdot \h{r}$ and $t' := \h{q} \cdot t \cdot q$ we find that 
$$
\xymatrix {
\D \ar@/_1.1pc/[rr]_{t'}^{\top} \ar@/^1.1pc/[rr]^{\y'}_{\top} & & \h{\D} \ar[ll]|{c'}
}
$$
so $\D$ is totally distributive.

If $t$ and $q$ are cartesian, then since $\h{q}$ is also cartesian, $t' = \h{q} \cdot t \cdot q$ is cartesian and hence $\D$ is lex totally distributive.
\end{proof}

\section{A construction of Johnstone-Joyal} \label{sec:const_of_jj}

Let $\X$ be an ind-small continuous category, and let $\A$ be a small ind-dense subcategory of $\X$.  We now recall from \cite{JJ} an explicit manner of constructing a quasi-injective topos $\F$ such that $\X$ is equivalent to the category of points of $\F$.

Firstly, there is an associated functor $W:\X^\op \times \X \rightarrow \Set$, given by
$$W(x,y) := \Ind{\X}(\m x,wy),\;\;\;\;x,y \in \X\;.$$
The elements of $W(x,y)$ are called \textit{wavy arrows} from $x$ to $y$ in $\X$.  Johnstone and Joyal \cite{JJ} show that this functor $W$, when viewed as a profunctor $W:\X \mrelarrow \X$, underlies an \textit{idempotent profunctor comonad} on $\X$, and that the restriction $V:\A^\op \times \A \rightarrow \Set$ of $W$ is again an idempotent profunctor comonad on $\A$.  In the latter case, since $\A$ is small, this means precisely that $V:\A \mrelarrow \A$ is an idempotent comonad on $\A$ in the bicategory $\Prof$ of small categories, profunctors, and morphisms of profunctors.  Further, $V$ is \textit{left-flat}, meaning that for each $y \in \A$, $V(-,y):\A^\op \rightarrow \Set$ is a flat presheaf.

Recall that for small categories $\C,\D$, each profunctor $M:\C \mrelarrow \D$ (by which we mean a functor $M:\C^\op \times \D \rightarrow \Set$) gives rise to a cocontinuous functor $\widetilde{M}:[\C,\Set] \rightarrow [\D,\Set]$.  Indeed, $\widetilde{M}$ is the left Kan extension along the Yoneda embedding $\C^\op \rightarrow [\C,\Set]$ of the transpose $\C^\op \rightarrow [\D,\Set]$ of $M$.  This passage defines an equivalence of the bicategory $\Prof$ with another bicategory, in fact a 2-category, whose objects are again all small categories, but whose 1-cells $\C \rightarrow \D$ are all cocontinuous functors $[\C,\Set] \rightarrow [\D,\Set]$, and whose 2-cells are all natural transformations.

Hence our idempotent comonad $V:\A \mrelarrow \A$ in $\Prof$ determines an idempotent comonad $\widetilde{V}:[\A,\Set] \rightarrow [\A,\Set]$.  Moreover, since $V(-,y):\A^\op \rightarrow \Set$ is flat for each $y \in \A$, it follows that $\widetilde{V}$ preserves finite limits and so is said to be a \textit{cartesian comonad}.  Further, since $\widetilde{V}$ is also cocontinuous, $\widetilde{V}$ is the inverse-image part of a geometric morphism:

\begin{DefSub}
Given an ind-small continuous category $\X$ with a small ind-dense subcategory $\A$, the \textit{associated geometric endomorphism} is defined to be the geometric morphism $m_{\A,\X}:[\A,\Set] \rightarrow [\A,\Set]$ whose inverse-image part is \textit{the associated idempotent comonad} $m_{\A,\X}^* =\widetilde{V}$.
\end{DefSub}

\begin{PropSub} \label{thm:jjassoctopos}
\emph{(Johnstone-Joyal \cite{JJ}).}
Let $\X$ be an ind-small continuous category, and let $\A$ be a small ind-dense subcategory of $\X$.  Let $[\A,\Set] \rightarrow \F \rightarrow [\A,\Set]$ be a factorization of the associated geometric endomorphism $m_{\A,\X}$ into a geometric surjection followed by a geometric inclusion.  Then $\F$ is a quasi-injective topos whose category of points of is equivalent to $\X$.  Further, if $\X$ is cocomplete, then we may take $\A$ to be finitely cocomplete, and it follows that $\F$ is an injective topos.
\end{PropSub}

\section{Totally distributive toposes from continuous categories}

We now show that the toposes corresponding to continuous categories under the equivalence of Theorem \ref{thm:injcts} are totally distributive, so that every quasi-injective topos is totally distributive.

\begin{LemSub} \label{thm:left_adjoint_idempotent_comonad}
Let $i:\C \rightarrow \D$ be a fully faithful functor with a right adjoint $r$, and suppose that the induced comonad $i \cdot r$ on $\D$ has a right adjoint $n$.  Then $r$ has a right adjoint $s := n \cdot i$, so that $i \dashv r \dashv s$.
\end{LemSub}
\begin{proof}
$$\C(r(d),c) \cong \D(i \cdot r(d),i(c)) \cong \D(d,n \cdot i(c)) = \D(d,s(c))\;,$$
naturally in $d \in \D, c \in \C$.
\end{proof}

\begin{LemSub} \label{thm:assoctopos}
Let $\X$ be an ind-small continuous category, let $\A$ be a small ind-dense subcategory of $\X$, and let $i:\F \hookrightarrow [\A,\Set]$ be the coreflective embedding induced by the associated idempotent comonad $m_{\A,\X}^*$ on $[\A,\Set]$ (so that $\F$ is the category of fixed points of $m_{\A,\X}^*$).  Then 
\begin{enumerate}
\item i preserves finite limits;
\item The right adjoint $r:[\A,\Set] \rightarrow \F$ to $i$ has a further right adjoint $s$, so that 
$$
\xymatrix {
\F \ar@/_1.1pc/[rr]_i^{\top} \ar@/^1.1pc/[rr]^s_{\top} & & {[\A,\Set]\;;} \ar[ll]|r
}
$$
\item $\F$ is a quasi-injective topos whose category of points is equivalent to $\X$;
\item If $\X$ is cocomplete, we may take $\A$ to be finitely cocomplete, and $\F$ is then an injective topos.
\end{enumerate}
\end{LemSub}
\begin{proof}
Since $\F$ is isomorphic to the category of coalgebras of the cartesian comonad $m_{\A,\X}^*$, $\F$ is an elementary topos, and the forgetful functor $i:\F \hookrightarrow [\A,\Set]$ is the inverse-image part of a geometric surjection $p:[\A,\Set] \twoheadrightarrow \F$; see, e.g., \cite{Jele}, A4.2.2.  Further, the idempotent comonad $i \cdot r = m_{\A,\X}^*$ has a right adjoint ${m_{\A,\X}}_*$, so we deduce by Lemma \ref{thm:left_adjoint_idempotent_comonad} that $r$ has a right adjoint $s$, so that $i \dashv r \dashv s$.  In particular, $r$ is left adjoint and cartesian, so we obtain a geometric morphism $q:\F \rightarrow [\A,\Set]$ with $q^* = r$ and $q_* = s$.  Since $i \dashv r \dashv s$ and $i$ is fully faithful, it follows that $s = q_*$ is also fully faithful, so $q:\F \rightarrow [\A,\Set]$ is a geometric inclusion.  Further, the composite $[\A,\Set] \xrightarrow{p} \F \xrightarrow{q} [\A,\Set]$ is $m_{\A,\X}$, or, more precisely, has inverse-image part $(q \cdot p)^* = p^* \cdot q^* = i \cdot r = m_{\A,\X}^*$.  Hence 3 and 4 follow from Proposition \ref{thm:jjassoctopos}.
\end{proof}

\begin{DefSub}
For an ind-small continuous category $\X$ and a small ind-dense subcategory $\A$ of $\X$, we call the topos $\F$ of Lemma \ref{thm:assoctopos} the \textit{associated topos}.
\end{DefSub}

\begin{LemSub} \label{thm:assoctopostd}
Let $\X$ be an ind-small continuous category, so that $\X$ has some small ind-dense subcategory $\A$.  Then the the associated topos $\F$ is totally distributive.  If $\X$ is also cocomplete, then we may take $\A$ to be finitely cocomplete, and it follows that $\F$ is lex totally distributive.
\end{LemSub}
\begin{proof}
By Lemma \ref{thm:assoctopos}, we have adjunctions
$$
\xymatrix {
\F \ar@/_1.1pc/[rr]_i^{\top} \ar@/^1.1pc/[rr]^s_{\top} & & [\A,\Set] \ar[ll]|r
}
$$
with $i,s$ fully faithful and $i$ cartesian.  By Proposition \ref{thm:preshtd}, $[\A,\Set]$ is totally distributive, so we deduce by Lemma \ref{thm:mrwlemma} that $\F$ is totally distributive.  If $\X$ is also cocomplete, then we can take $\A$ to be finitely cocomplete, so $\A^\op$ is finitely complete and hence, by \ref{thm:preshtd}, $\h{\A^\op} = [\A,\Set]$ is lex totally distributive, so we deduce by \ref{thm:mrwlemma} that $\F$ is lex totally distributive.
\end{proof}

\begin{ThmSub} \label{thm:injimpliesltd}
Every quasi-injective topos is totally distributive, and every injective topos is lex totally distributive.
\end{ThmSub}
\begin{proof}
Given a quasi-injective topos $\E$, Theorem \ref{thm:injcts} entails that the category of points $\X := \pt(\E)$ of $\E$ is an ind-small continuous category.  Taking any small ind-dense subcategory $\A$ of $\X$, the associated topos $\F$ is a quasi-injective topos whose category of points is equivalent to $\X$, so by Theorem \ref{thm:injcts} we deduce that $\E$ is equivalent to $\F$.  But the latter topos is totally distributive by Lemma \ref{thm:assoctopostd}, and total distributivity is clearly invariant under equivalences, so $\E$ is totally distributive.  The second statement may be deduced analogously.
\end{proof}

\section{Totally distributive categories as essential localizations}

\begin{PropSub} \label{thm:resttogen}
Let $\E$ be a totally distributive category with a small dense generator $i:\G \hookrightarrow \E$.  We then conclude the following: 
\begin{enumerate}
\item There are adjunctions
$$
\xymatrix {
\E \ar@/_1.1pc/[rr]_{t'}^{\top} \ar@/^1.1pc/[rr]^{\y'}_{\top} & & {\h{\G}} \ar[ll]|{c'}
}
$$
with $\y'$ and $t'$ fully faithful, where $\y'$ is the composite $\E \xrightarrow{\y} \h{\E} \xrightarrow{\h{i}} \h{\G}$.
\item $\E$ is an essential subtopos of $\h{\G}$ and, in particular, a Grothendieck topos.
\item If $\E$ is lex totally distributive, then $t':\E \rightarrow \h{\G}$ preserves finite limits.
\end{enumerate}
\end{PropSub}
\begin{proof}
We let
$$c' := c \cdot \forall_i = (\h{\G} \xrightarrow{\forall_i} \h{\E} \xrightarrow{c} \E)\;,$$
$$t' := \h{i} \cdot t = (\E \xrightarrow{t} \h{\E} \xrightarrow{\h{i}} \h{\G})\;,$$
where $\forall_i:\h{\G} \rightarrow \h{\E}$ is the functor given by right Kan extension along $i^\op:\G^\op \hookrightarrow \E^\op$.  Since $\h{i} \dashv \forall_i$ and $t \dashv c$, we have that $t' = \h{i} \cdot t \dashv c \cdot \forall_i = c'$.  Since $i:\G \hookrightarrow \E$ is fully faithful, the counit of the adjunction $\h{i} \dashv \forall_i$ is an isomorphism (e.g., by \cite{Ma}, X.3.3), so $\forall_i$ is fully faithful.

Observe that the diagram
$$
\xymatrix {
\E     \ar[d]_\y \ar[r]^\y  & \h{\E} \\
\h{\E} \ar[r]_{\h{i}}       & \h{\G} \ar[u]_{\forall_i}
}
$$
commutes up to isomorphism, since the density of $\G$ in $\E$ gives us exactly that
$u \cong \int^{g \in \G} \E(g,u) \cdot g$
naturally in $u \in \E$, so
$$(\y v)u = \E(u,v) \cong \E(\int^{g \in \G} \E(g,u) \cdot g,v) \cong \int_{g \in \G} [\E(g,u),\E(g,v)] = ((\forall_i \cdot \h{i} \cdot \y)v)u$$
naturally in $u,v \in \E$.  

We find that $c' = c \cdot \forall_i \dashv \h{i} \cdot \y = \y'$, since by using the adjointness $c \dashv \y$, the commutativity of the above diagram, and the fact that $\forall_i$ is fully faithful, we deduce that
$$\E(c \cdot \forall_i(G),v) \cong \h{\E}(\forall_i(G),\y v) \cong \h{\E}(\forall_i(G),\forall_i \cdot \h{i} \cdot \y(v)) \cong \h{\G}(G,\h{i} \cdot \y(v))$$
naturally in $G \in \h{\G}, v \in \E$.  

Since $\G$ is a dense generator for $\E$ we have that $\y'$ is fully faithful, and since $t' \dashv c' \dashv \y'$ it follows that $t'$ is fully faithful as well.

If $\E$ is lex totally distributive, then $t$ preserves finite limits, so since $\h{i}$ preserves all limits, $t' = \h{i} \cdot t$ preserves finite limits.
\end{proof}

\begin{ThmSub} \label{thm:lextdinjective}
Let $\E$ be a lex totally distributive category with a small dense generator.  Then $\E$ is an injective Grothendieck topos.
\end{ThmSub}
\begin{proof}
By \ref{thm:resttogen} we know that $\E$ is a Grothendieck topos, and it follows from Giraud's Theorem that there exists a \textit{finitely complete} small dense full subcategory $\G$ of $\E$.  (Indeed, this follows readily from 4.1 and 4.2 in the Appendix of \cite{MM}, for example).  We have adjunctions $t' \dashv c' \dashv \y'$ as in Proposition \ref{thm:resttogen}, with $\y'$ fully faithful and $t'$ cartesian.  Hence we obtain geometric morphisms $s:\E \rightarrow \h{\G}$ and $r:\h{\G} \rightarrow \E$ with
$s_* = \y'$, $s^* = c'$, $r_* = c'$, $r^* = t'$, since $c'$ is right adjoint and hence cartesian.  Further, since $\y'$ is fully faithful and $c' \dashv \y'$, we have that
$$(r \cdot s)_* = r_* \cdot s_* = c' \cdot \y' \cong 1_\E \;,$$
so $\E$ is a (pseudo-)retract of the presheaf topos $\h{\G}$ by geometric morphisms, and the result follows by \ref{thm:injtop_as_retracts}.
\end{proof}

Hence Theorem \ref{thm:lex_td_sdg_equals_inj_top} is proved.  To prove Theorem \ref{thm:charn_of_td_sdg}, it remains only to show the following:

\begin{PropSub}
Let $\E$ be an essential subtopos of a presheaf topos $\h{\C}$ (with $\C$ small).  Then $\E$ is totally distributive and has a small dense generator.
\end{PropSub}
\begin{proof}
There is a geometric inclusion $s:\E \rightarrow \h{\C}$ whose inverse-image functor $s^*:\h{\C} \rightarrow \E$ has a left adjoint $s_!$.  Hence we have $s_! \dashv s^* \dashv s_*$ with $s_!$ and $s_*$ fully faithful, so $\E$ is totally distributive, by Lemma \ref{thm:mrwlemma}.
\end{proof}

\bibliographystyle{amsplain}
\bibliography{td1}

\end{document}